\documentclass[12pt]{amsart}
\usepackage{amssymb}

\newtheorem{Thm}{Theorem}[section]

\newtheorem{Cor}[Thm]{Corollary}
\newtheorem{Lem}[Thm]{Lemma}
\newtheorem{Prop}[Thm]{Proposition}

\def\ldots{\mathinner{\ldotp\ldotp\ldotp}}
\def\ldots{\mathinner{\cdotp\cdotp\cdotp}}

\def \R{\rm {\bf R}}

\def\Ave{\mathop{\text {Ave}}}
\def\r{\right}

\begin{document}

\title{Quotients of finite-dimensional quasi-normed spaces}
\author{N. J. Kalton}
\address{Department of Mathematics \\
University of Missouri-Columbia \\
Columbia, MO 65211
}

\email{nigel@math.missouri.edu}

\author{A. E. Litvak}
\address{Department of Mathematics \\
Technion, Haifa\\
Israel,  32000\\
}

\email{alex@math.technion.ac.il}

\thanks{The first author was supported by NSF grant DMS-9870027 and the
second
   author was supported by a Lady Davis Fellowship}


\begin{abstract} We study the existence of cubic quotients of
finite-dimensional quasi-normed spaces, that is, quotients well
isomorphic to
$\ell_{\infty}^k$ for some $k.$  We give two results of this nature.  The
first guarantees a proportional dimensional cubic quotient when the
envelope is cubic; the second gives an estimate for the size of a cubic
quotient in terms of a measure of non-convexity of the quasi-norm.
\end{abstract}

\maketitle

\section{Introduction}
It is by now well-established that many of the core results in the local
theory of Banach spaces extend to quasi-normed spaces (cf. \
\cite{BBP1}, \cite{BBP2}, \cite{D}, \cite{GK}, \cite{GL}, \cite{K1},
\cite{K2}, \cite{KT}, \cite{L1}, \cite{LMP}, \cite{M} for example).
 In this note we give two results on the local theory of quasi-normed
spaces which are of interest only in the non-convex situation.

Let us introduce some notation.  Let $X$ be a real finite-dimensional
vector space. Then a
$p$-norm $\|\cdot\|$ on $X$, $p\in (0,1]$,
 is a map $x\mapsto \|x\|\ (X\mapsto \mathbb R)$
so that: \newline
(i) $\|x\|>0$ if and only if $x\neq 0.$\newline
(ii) $\|\alpha x\|=|\alpha|\|x\|$ for $\alpha\in\mathbb R$ and $x\in X$.
\newline
(iii) $\|x_1+x_2\|^p\le \|x_1\|^p+\|x_2\|^p$ for $x_1,x_2\in X.$
\newline
Then $(X,\|\cdot\|)$ is called a $p$-normed space.  For the purposes of this
paper a quasi-normed space is always assumed to be a $p$-normed space for
some $p\in (0, 1]$ (note that by Aoki-Rolewicz theorem on quasi-normed space
one can introduce an equivalent $p$-norm (\cite{KPR}, \cite{Ko}, \cite{R})).
The set $B_X=\{x:\|x\|\le 1\}$ is the unit ball of $X$.
The closed convex hull of $B_X$, denoted by $\hat B _X$,
 is the unit ball of a norm $\|\cdot\|_{\hat X}$ on $X$; the corresponding normed
space, $\hat X$, is called the Banach envelope of $X$.

The set $B$ is called $p$-convex if for every $x$, $y \in B$ and every positive
$\lambda$, $\mu$ satisfying $\lambda ^p + \mu ^p = 1$ one has
$\lambda x + \mu y \in B$. Clearly, the unit ball of $p$-normed space is a
$p$-convex set and, vise versa, a closed centrally-symmetric $p$-convex set is
the  unit ball of some $p$-norm  provided that it is bounded and $0$
belongs to its interior.

If $X$ and $Y$ are $p$-normed spaces (for some $p$) then the Banach-Mazur
distance $d(X,Y)$ is defined as $\inf\{\|T\|\|T^{-1}\|  \}$, where the infimum is taken
over all linear isomorphisms $ T:X\to Y$.
  We let $d_{B_X} = d_X = d(X,\ell_2^{\dim X})$
and $\delta_{B_X} = \delta_X = d(X,\hat X).$
It is clear that $\delta_X$ is measure of
non-convexity; in fact $\delta_X=\inf\{d(X,Y):\ Y \text{ is a Banach
space}\}.$

We now describe our main results.  In Section \ref{cube} we investigate
quasi-normed spaces $X$ such that $\hat X$ satisfies an estimate
$d(\hat X,\ell_{\infty}^{\dim X})\le C.$  It has been known for some time that
non-trivial examples of this phenomenon exist \cite{K}.  In geometrical
terms this means that the convex hull of the unit ball of $X$ is close to
a cube.  We show using combinatorial results of Alesker, Szarek and
Talagrand \cite{A}, \cite{ST} based on the Sauer-Shelah Theorem
\cite{Sa}, \cite{Sh} that $X$ then has a proportional dimensional
quotient $E$ satisfying an estimate $d(E,\ell_{\infty}^{\dim E})\le C'.$
 A much more precise statement is given in Theorem \ref{pnormedmain}.
We then use this result in Section~\ref{cube2} to prove that a $p$-normed
space $X$ has a quotient $E$ with $\dim E\ge c_p \ln \delta_X /
(\ln \ln \delta_X )$ and  $d(E,\ell_{\infty}^{\dim E})\le C_p$  where
$0<c_p ,C_p<\infty$ are constants  depending on $p$ only.
Again a more precise statement is given in Theorem~\ref{final}.

In developing these results, we  found it helpful to use the notion of a
geometric hull of a subset of $\mathbb R^n.$  Thus instead of considering
a $p$-convex set $B_X$ we consider an arbitrary compact spanning set $S$
and then compare the absolutely convex hull $\Delta S$ with certain
subsets $\Gamma_{\theta}S$ which can be obtained from $S$ by
geometrically converging series.  Precisely $x\in\Gamma_{\theta}S$,
$\theta \in (0,1)$, if and only if
$x=(1-\theta)\sum_{n=0}^{\infty}\theta^n\lambda_ns_n$ where
$s_n\in S$ and $|\lambda_n|\le 1.$ Note  that $\Gamma _{\alpha} \subset
\frac{1-\alpha}{1-\theta} \Gamma _{\theta} $ for every $0 < \alpha< \theta < 1$.
 Our results can be stated in terms of
estimates for the speed of convergence of $\Gamma_{\theta}S$ to $\Delta
S$ as $\theta\to 1.$  In this way we can derive results which are
independent of $0<p<1$ and then obtain results about $p$-normed spaces as
simple Corollaries.  We develop the idea of the geometric hull in Section
\ref{approx} and illustrate it by restating the quotient form of
Dvoretzky's theorem in this language.

\section{Approximation of convex sets} \label{approx}

Let $S$ be a subset of $\mathbb R^n.$  Denote by $\Delta S$ the
absolutely convex hull of $S$ and by $\tilde S$
 the star-shaped hull of $S,$ i.e.
$\tilde S=\{\lambda x:\ |\lambda|\le 1,\ x\in S\}.$
  For each $m\in\mathbb N$ we define
$\Delta_mS$ to be the set of all vectors of the form
$\frac1m(\lambda_1x_1+\cdots+\lambda_mx_m) $ where $|\lambda_k|\le 1$ and
$x_k\in S$ for $1\le k\le m.$  If $0<\theta<1$ we define  {\it the
$\theta$-geometric hull} of $S$,
$\Gamma_{\theta}S$ to be the set of all vectors of the form
$(1-\theta)\sum_{k=0}^{\infty}\lambda_kx_k$ where $|\lambda_k|\le
\theta^k$ and $x_k\in S$ for $k=0,1,\cdots.$

\begin{Lem}\label{pconv}  Let $S$ be a $p$-convex closed set where $0<p<1.$
Then  for $0<\theta<1$ we have $$\Gamma_{\theta}S \subset
\left( p^{-1/p} (1-\theta)^{1-1/p} \right) S.$$
\end{Lem}

\begin{proof}  This follows easily from:
\begin{equation*}
 \frac{1-\theta}{(1-\theta^p)^{1/p}} \le
p^{-1/p}(1-\theta)^{1- 1/p}\end{equation*}
which in turns from the estimate $$ \theta^p\le
1-p(1-\theta).$$\end{proof}

\begin{Lem}\label{approx2} If $\frac13< \theta< 1$ and $m\in\mathbb
N$ then  $$ \Gamma_{\theta}\Delta_mS \subset
\frac{2\theta}{3\theta-1}\Gamma_{\theta^{\frac1m}}S.$$\end{Lem}

\begin{proof} Note that $$\Delta_mS \subset \frac1m
\theta^{\frac1m-1}\sum_{k=0}^{m-1}\theta^{\frac{k}{m}}\tilde S . $$
 Hence $$ \Gamma_{\theta}\Delta_mS \subset
\frac{1-\theta}{m(1-\theta^{\frac1m})}\theta^{\frac1m-1}
\Gamma_{\theta^{\frac1m}}S.$$
Now observe
\begin{equation}\begin{align*}
\frac{1-\theta}{m(1-\theta^{\frac1m})}\theta^{\frac1m-1}&=
\frac{\theta^{-1}-1}{m(\theta^{-\frac1m}-1)}\\
&\le \frac{\theta^{-1}-1}{|\ln \theta|}\\
&\le \frac{2}{3-\theta^{-1}}.\end{align*}\end{equation}
This completes the proof.\end{proof}

 In this section, we make a few simple observations on the geometric
hulls $\Gamma_{\theta}S.$  Let us suppose that $S$ is compact and
spanning so that $\Delta S$ coincides with the unit ball $B_X$ of a
Banach space $X,\ \|\cdot\|_X.$
 Given $q \in [1, 2]$ let $T_q = T_q (X)$ denote the equal-norm type $q$
constant, i.e.  the smallest constant satisfying
$$ \Ave_{\epsilon_k=\pm1}\left\|\sum_{k=1}^N\epsilon_kx_k\r\|_X\le
T_q N^{1/q} \max_{1\le k\le N}\|x_k\|$$ for every $N$.
Given an integer $N$   let $b_N$ denote the least constant so that
$$ \inf_{\epsilon_k=\pm1}\left\|\sum_{k=1}^N\epsilon_kx_k\r\|_X\le
b_NN\max_{1\le k\le N}\|x_k\|.$$
 Given a set $A$ by $|A|$ we denote the cardinality of $A$.

The following Lemma abstracts the idea of \cite{GK}, Lemma 2.

\begin{Lem}\label{type1} Suppose $\frac13<\theta<1$, and let $m=m(S)$ be
an integer such that $\sum_{k=1}^{\infty}b_{2^km}\le \theta.$
Then
$$ \Delta S \subset
\frac{2\theta}{(3\theta-1)(1-\theta)}\Gamma_{\theta^{\frac1m}}S.$$
\end{Lem}

\begin{proof} Suppose $N\in\mathbb N$ and suppose $u\in \Delta_{2N}S.$
Then $u=\frac{1}{2N}(x_1+\cdots +x_{2N})$ where $x_k\in \tilde S.$
Hence there is a choice of signs $\epsilon_k=\pm1$ with
$|\{\epsilon_k= - 1\}|\le N$ and
$$ \left\|\sum_{k=1}^{2N} \epsilon_k x_k\right\|_X \le 2Nb_{2N}.$$
Let $v=\frac1N(\sum_{\epsilon_k=1}x_k).$  Then
$ \|u-v\|_X \le b_{2N}$.  Hence $\Delta_{2N}S \subset
\Delta_NS+b_{2N}\Delta S.$  Iterating we get
$$ \Delta_{2^km} S \subset \Delta_mS +\sum_{j=1}^kb_{2^jm}\Delta S$$
which leads to
$$ \Delta S  \subset \Delta_mS +\theta\Delta S$$ which implies
$$ \Delta S \subset (1-\theta)^{-1}\Gamma_{\theta}\Delta_mS \subset
\frac{2\theta}{(3\theta-1)(1-\theta)}\Gamma_{\theta^{\frac1m}}S.$$\end{proof}

\begin{Prop}\label{type2}
(i) Suppose $1<q \le 2$ and $q'$ be such that $1/q + 1/q' =1$. Then for
$$
\theta = 1 - \frac{1}{4} \left( \frac{2^{1/q'} -1}{2 T_q} \right) ^{q'}
$$ we have
$\Delta S\subset 12\Gamma_{\theta}S.$\newline
(ii) There exists  constant $C<\infty$ so that if $m$ is the
largest integer such that $X$ has a subspace $Y$ of dimension
$m$ with $d(Y,\ell_1^m)\le 2$ then
$\Delta S \subset 8 \Gamma_{\theta}S$ for
$$\theta =1-\frac12 \left(C m \r) ^{ - C \log\log (Cm) }.$$
\end{Prop}

\medskip

\noindent
{\bf Remark.}
We conjecture that the sharp estimate in (ii) is $\theta  = 1-c / m.$

\begin{proof}(i)
 Observe that $b_N\le T_q N^{\frac1q-1}$. Hence
$$ \sum_{k=1}^{\infty}b_{2^kN} \le T_q
N^{-\frac{1}{q'}}(2^{\frac{1}{q'}}-1)^{-1}.$$
Let $N$ be the largest integer so that the right-hand side is
at most $\frac12$. Applying Lemma \ref{type1} with $\theta _0 =1/2$
we obtain
$$\Delta S\subset 4 \Gamma_{2^{-1/N}} S.$$
The result follows, since
$$
 \frac{1}{N} \leq \left( \frac{2^{1/q'} -1}{2 T_q} \right) ^{q'} \leq  \frac{1}{N-1}
\ \ \mbox{ and } \ \ \Gamma _{\alpha} \subset
\frac{1-\alpha}{1-\theta} \Gamma _{\theta}
$$
for $\alpha< \theta$.

In (ii) we note first by a result of Elton \cite{E} (see also \cite{T}
for a sharper version) there exist universal constants
$1/2 \leq c_0 <1$ and $C\geq 1$ so that $b_{N_0}< c_0$
for some $N_0\le Cm.$

Recall simple properties of the numbers $b_k$. Clearly, for every
$k$, $l$ one has $b_{k l} \leq b_k b_l $ and $(k+l) b_{k+l} \leq
k b_k + l b_l$. Thus if  $b_k \leq c_0 <1$ then $b_l \leq c =(1+c_0)/2 <1$
for every $k \leq l \leq 2k$.
Therefore we may
 suppose that $N_0$ is a power of two, say $N_0 = 2^q$, $q\geq 1$, and
$b_{ N_0  } \leq c <1$. Since $b_l \leq 1$ for every $l$, we get
 $b_{ N_0^s l } \leq c^s$ for every integers $s\geq 1$, $l\geq 0$.
Then, taking $N= N_0^r$ for some $r\geq 1$ we have
$$
 \sum_{k=1}^{\infty}b_{2^k N} = \sum_{j=0}^{\infty} \sum_{l=1}^{rq}
 b_{2^{rq+jrq+l}} \leq r q   \sum_{j=1}^{\infty} c^{j r}
\leq 2 r q c^r \leq 1/2
$$
provided $r \geq c_1 \ln q$ with appropriate absolute constant $c_1$.

Now take $r$ to be smallest integer larger than
$ c_1 \ln q = c_1 \ln \log_2 N_0$. Then by Lemma~\ref{type1} we obtain
$$
   \Delta S \subset 4  \Gamma _{2^{-1/N}} S
$$
for $N\sim  (C'm)^{C' \log\log (C'm)} $ and the result follows.
\end{proof}

\begin{Cor}
 There are absolute constants $c$, $C>0$ so that  if $X$ is a $p$-normed
space then there exists a subspace $Y$ in the envelope $\hat X$
such that dimension of $Y$ is
$$ m \geq c p \exp \left\{ \frac{\ln A}{ \ln \ln A}\r\} ,$$
where $A= C (\delta _X )^{p/(1-p)}$,
and
$$ d\left(Y , \ell _1^m\r) \leq 2 .$$
\end{Cor}

\begin{proof}
 Let $S=B_X$ and let $m$ be as in Proposition~\ref{type2}.
Then by the proposition we
have $\Delta B_X \subset 8 \Gamma _{\theta} B_X$ with
$$ \theta =1-\frac12 \left(C m \r) ^{ - C \log\log (Cm) }.    $$
 Thus by Lemma~\ref{pconv} we obtain
$$  \Delta B_X  \subset 8 p^{-1/p} 2^{-1+1/p} \left(C m \r) ^{ - \left(1-1/p \r)
C \log\log (Cm) } B_X,$$
i.e.
$$ \delta _X \leq \left( C' m / p \r)^{-\left(1-1/p \r) C \log\log (Cm) } .$$
That implies the result.
\end{proof}

Let us conclude this section with a very simple form of Dvoretzky's
theorem recast in this language:

\begin{Thm}\label{Dvor} Let $\eta < 1/3$.
 There is an absolute  constant $c>0$ so that if $S$ is a
compact spanning subset of $\mathbb R^n$ then there is a projection $P$
of rank at least $c \eta^2 \log n$ such that
$$d_{\Gamma_{\theta}PS}\le \frac{1+\eta}{1-\theta} $$
for every  $  \sqrt{3 \eta}\leq \theta < 1.$
\end{Thm}

\medskip

\noindent
{\bf Remark 1.} Let $\epsilon \leq 6/7$.
Setting $\theta = \sqrt{3 \eta} =  \epsilon /2 $ we observe that
there is an absolute  constant $c>0$ so that if $S$ is a
compact spanning subset of $\mathbb R^n$ then there is a projection $P$
of rank at least $c \epsilon^4 \log n$ such that
$$d_{\Gamma_{ \epsilon /2}PS}\le 1+\epsilon. $$

\noindent
{\bf Remark 2.}
The ``quotient form'' of Dvoretzky's theorem for quasi-normed spaces
is essentially known and follows very easily from results in \cite{GK}
(see e.g. \cite{GL} for the details).

\begin{proof}
 By the sharp form of Dvoretzky's Theorem
(Theorem 2.9 in \cite{G}) there is a projection $P$ of rank at least
$c\eta^2\log n$ so that $d_{\Delta(PS)}\le 1+\eta.$
Let $Y=P\mathbb R^n$ and introduce an inner-product norm
$\|\cdot\|$ on $Y$ so that
${\mathcal E}\subset \Delta(PS)\subset(1+\eta) {\mathcal E}$
where ${\mathcal E}=\{y:\ (y,y)\le 1\}.$  If $y\in \mathcal E$ with
$\|y\|=1$ there exists $u\in PS \cup (-PS)$ with $(y,u)\ge 1.$
Since $\|u\|\le 1+\eta$ we obtain
$\|y-u\|\le (2\eta+\eta^2)^{1/2}\le \sqrt{3 \eta }.$
Hence
$$ {\mathcal E}\subset PS\cup (-PS)+\sqrt{3 \eta }\ {\mathcal E}$$
which implies, for any $\theta \geq \sqrt{3 \eta }$,
$$ (1-\theta ) {\mathcal E}\subset \Gamma_{\theta}PS
\subset (1+\eta) {\mathcal E}.$$
Hence
$$ d_{\Gamma_{\theta } PS} \le \frac{1+\eta}{1-\theta },$$
which proves the theorem.
\end{proof}

\section{Approximating the cube}\label{cube}

Let $n$ be an integer. By $[n]$ we denote the set $\{1, ..., n\}$.
The $n$-dimensional cube we denote by $B^{\infty} = B^{\infty}_n$.
$D_n$ denotes the extreme points of the cube, i.e. the set $\{1, -1\} ^n$.
Given a set $\sigma \subset [n]$ by $P_{\sigma}$ we denote the
coordinate projection of $\mathbb R^n$ onto $\mathbb R^{\sigma}$, and
we denote  $B^{\infty}_{\sigma}:= P_{\sigma} B^{\infty}_n$,
$D_{\sigma} := P_{\sigma} D_n$.
As above $|A|$ denotes the cardinality of a set $A$. As usual $\| \cdot \| _2$
and $\| \cdot \| _{\infty}$ denote the norm in $\ell _2$ and $\ell _{\infty}$
correspondingly.

\begin{Thm}\label{alesker} There are constants $c>0$ and $0<C<\infty$ so
that for every
$\epsilon>0$, if $S\subset D_n$ with $|S|\ge 2^{n(1-c\epsilon)}$ then
there is a subset $\sigma$ of $[n]$ with $|\sigma|\ge (1-\epsilon ) n$ so
that
$$ D_{\sigma}\subset C\epsilon^{-1}P_{\sigma}(\Delta_N S)$$ for some
$N\le C\epsilon^{-2}.$\end{Thm}

\begin{proof} We will follow Alesker's argument in \cite{A}, which
is itself a refinement of Szarek-Talagrand \cite{ST}.  Alesker shows
that for a suitable choice of $c$, if $\epsilon=2^{-s}$
then one can find an increasing sequence of subsets $(\sigma_k)_{k=0}^s$
so that $P_{\sigma_0}(S)=D_{\sigma_0},$
$|\sigma_s|\ge
(1-2\epsilon)n$ and if $\tau_k=\sigma_k\setminus \sigma_{k-1}$ for
$k=1,2,\ldots,s$ then there exists $\alpha\in D_n$ so that
$$ P_{\tau_k}(S\cap
P_{\sigma_{k-1}}^{-1}(P_{\sigma_{k-1}}\alpha))=D_{\tau_k}.$$
It follows that if $a\in D_{\tau_k}$ there exists $x\in \Delta_2S$ with
$P_{\sigma_{k-1}}(x)=0$ and $P_{\tau_k}(x)=a.$

We now argue by induction that $D_{\sigma_k}\subset
a_k P_{\sigma_{k}}\Delta_{b_k}S$ where $a_k = 2^{k+1}-1$
and $b_k = 2^k a_k = 2\cdot 4^k - 2^k$. This
clearly holds if $k=0.$ Assume it
is true for
$k=j-1,$ where $1\le j\le s.$ Then if $a\in D_{\sigma_{j}}$ we can
observe that there exists
$x_1\in a_{j-1} \Delta_{b_{j-1}} S$ with
$P_{\sigma_{j-1}}x_1=P_{\sigma_{j-1}}a.$ Clearly,
$$P_{\tau_{j}}x_1\in a_{j-1} \Delta_{b_{j-1}} D_{\tau_j}.$$
Hence there exists $x_2\in a_{j-1} \Delta_{2 b_{j-1}}S$ with
$P_{\sigma_{j-1}}x_2=0$ and $P_{\tau_{j}}x_2=-P_{\tau_{j}} x_1.$
Finally pick $x_3\in \Delta_2S$ so that $P_{\sigma_{j-1}}(x_3)=0$ and
$P_{\tau_{j}}(x_3)=P_{\tau_{j}}a.$  Then
$P_{\sigma_{j}}(x_1+x_2+x_3)=a$ and
$$ x_1+x_2+x_3 \in a_{j-1} \Delta_{ b_{j-1}} S
+ a_{j-1} \Delta_{2 b_{j-1}} S +  \Delta_2 S $$
$$ \subset \frac{a_{j-1}}{2 b_{j-1}} \left( 4 b_{j-1} + 2^j \right)
 \Delta_{4 b_{j-1} + 2^j }S = a_{j} \Delta_{ b_{j}} S .$$
This  establishes the induction.

We finally conclude that
$D_{\sigma_s}\subset 2 (2^{s+1}-1) P_{\sigma_s} \Delta_{2 \cdot 4^s}S$
and this gives the result, as the case of general $\epsilon$ follows
easily.\end{proof}

\medskip

\noindent
{\bf Remark.} Slightly changing the proof one can show that
$ D_{\sigma}\subset C\epsilon^{-\alpha}P_{\sigma}(\Delta_N S)$
for $N\le C\epsilon^{-\alpha}$, where $\alpha = \log_2 3$.

\begin{Lem}\label{counting}  There exist absolute constants $c,C>0$ with
the following property.  Suppose
$0<\epsilon<1$ and $0<k<n$ are natural numbers with $k/n\ge
1-c\epsilon(1-\ln\epsilon )^{-1}$.  Let $S$ be a subset of $\mathbb R^n$
so that if $a\in D_n$ there exists $x\in S$ with $|\{i:x_i=a_i\}|\ge k.$
Then there is a subset $\sigma$ of $[n]$ with $|\sigma|\ge
(1-\epsilon)n$ and $D_{\sigma}\subset C\epsilon^{-1}\Delta_NP_{\sigma}S$
for some $N\le C\epsilon^{-2}.$ \end{Lem}

\begin{proof} Suppose $0<k<n$ and $1-k/n=t\epsilon (1-\ln\epsilon)^{-1}.$
We shall show that if $t$ is small enough we obtain the conclusion of the
lemma.  First pick a map $a\to\sigma(a)$ from $D_n\to 2^{[n]}$ so that
for each $a$, $|\sigma(a)|=k$ and there exists $x\in S$ with $x_i=a_i$
for $i\in\sigma(a).$  Then, by a simple counting argument we have the
existence of $\tau\in 2^{[n]}$ so that $|\tau|=k$ and if
$$
T=\{\alpha \in D_{\tau}:\ \exists a\in D_n,\ \sigma(a)=\tau,\
P_{\tau} a =\alpha\}$$ then $$|T|\ge \frac{2^n}{2^{n-k} \binom{n}{k}}.
$$
We can estimate
$$ \binom{n}{k}\le\left(\frac{n}{k}\right)^k\left(\frac{n}{n-k}\right)^{n-k}
 \le \left(\frac{ne}{n-k}\right)^{n-k}.$$
Hence for $t \leq 1/2$ we have
$$ \log_2\binom{n}{k} \le
 \frac{n  t  \epsilon }{\ln 2 \left( 1-\ln  \epsilon \r)  }
  \ln \left(\frac{e \ln (e/\epsilon )}{t \epsilon } \r)
\leq 3  k   t  \epsilon \left( 2 - \ln t \r) .$$
  It follows that
$$ |T|\ge 2^{k(1-C_t \epsilon)}, $$
where $C_t = 3  t \left( 2 - \ln t \r) $. Choosing $t$ such that $C_t \leq c/2$,
where $c$ is the constant from Theorem~\ref{alesker}, and applying this
theorem, we obtain the existence of $\sigma \subset \tau$,
$|\sigma | \geq (1-\epsilon /2) k \geq (1-\epsilon ) n$, with desired property.
\end{proof}

\begin{Thm}\label{main}
 There are  absolute constants $c,  C>0$ such that if
$\epsilon>0$ and $S$ is a subset of
$\mathbb R^n$ with $B^{\infty}\subset \Delta S\subset dB^{\infty}$ then
there is a subset $\sigma$ of $[n]$ with $|\sigma|\ge n(1-\epsilon)$ such
that
$$  B^{\infty}_{\sigma} \subset \left( C/\epsilon\r) \Gamma_{\theta}P_{\sigma}S$$
for $\theta = 1-cd^{-2}\epsilon^5(1-\ln \epsilon)^{-1}.$
\end{Thm}

\begin{proof} Let $\delta=c_1\epsilon$  and $m$ be the smallest integer
greater than $c_2d^2\epsilon^{-3}(1-\ln \epsilon)$, where
$c_1,c_2$ will be chosen later.

Suppose first that $a\in D_n .$  Then  we can find
$N\in\mathbb N$, $N\geq m$, and $x_1,\ldots,x_N\in S\cup (-S)$ so that
$$ \left\|a-\frac1N(x_1+\cdots+x_N)\r\|_{2}^2\le \frac{nd^2}{m}.$$
Let $\Omega$ be the space of all $m$-subsets of $[N]$ and let $\mu$ be
normalized counting (probability) measure on $\Omega.$  If
$(\xi_i)_{i=1}^N$ denote
the indicator functions $\xi(\omega)=1$ if $i\in\omega$ and $0$ otherwise
then
$$ \bold E(\xi_i)=\bold E(\xi_i^2)=\frac{m}{N}, \ \bold
E(\xi_i\xi_j)=\frac{m(m-1)}{N(N-1)}$$
if $i\neq j.$  Thus
$$ \bold E(\xi_i-\bold E(\xi_i))^2 = \frac{m}{N}-\frac{m^2}{N^2}$$
and
$$ \bold E((\xi_i-\bold E(\xi_i))(\xi_j-\bold
E(\xi_j))=\frac{m(m-1)}{N(N-1)}-\frac{m^2}{N^2}$$
if $i\neq j.$

Let $y=\frac1N(x_1+\cdots +x_N)$ so that
$y=\bold E(\frac1m\sum_{i=1}^N\xi_ix_i).$
 Then working in the $\ell_2$-norm we have
$$\bold
E\left(\left\|\frac1m\sum_{i=1}^N\xi_ix_i-y\r\|_2^2\r)=
\frac{N-m}{mN(N-1)}\sum_{i=1}^N\left\|x_i \r\|_2^2
-\frac{N-m}{m N^2(N-1)}\left\|\sum_{i=1}^Nx_i\r\|_2^2.$$
Hence
$$\bold
E\left(\left\|\frac1m\sum_{i=1}^N\xi_ix_i-y\r\|_2^2\r)\le \frac{nd^2}{m}.$$
Since $\|y-a\|_2^2 \le \frac{nd^2}{m}$ we have
$$\bold
E\left(\left\|\frac1m\sum_{i=1}^N\xi_ix_i-a\r\|_2^2\r)\le 4\frac{nd^2}{m}.$$

We now suppose that for each $\omega\in \Omega$ we have $|\{j:\
|\frac{1}{m} \sum_{i=1}^N\xi_ix_i(j)-a(j)|>\delta\}|>4d^2n/(m\delta^2).$
Then we get an immediate contradiction.  We conclude that for each
$a\in D_n$ there exists $x_a\in \Delta_mS$ such that $|x_a(j)-a(j)|\le \delta$
for at least $n(1-2c_1^{-2}c_2^{-1}\epsilon(1-\log\epsilon)^{-1})$ choices
of $j.$  Let $y_a(j)=a(j)$ if $|x_a(j)-a(j)|\le \delta$ and $y_a(j)=x_a(j)$
otherwise so that $\|y_a-x_a\|_{\infty}\le \delta.$

Now suppose $c_2$ is chosen as a function of $c_1$ so that we can apply
Lemma \ref{counting} to obtain the existence of a set $\sigma\subset
[n]$ with
$|\sigma|\ge n(1-\epsilon)$ and so that
$$D_{\sigma}\subset C\epsilon^{-1}P_{\sigma}\Delta_N\{y_a:a\in D_n\}$$
where $C$ is an absolute constant, and $N\le C\epsilon^{-2}$.  Then
$$ D_{\sigma} \subset C\epsilon^{-1}P_{\sigma}\Delta_{Nm}S +
C\epsilon^{-1}\delta B^{\infty}_{\sigma}.$$
Recall that $C\epsilon^{-1}\delta =Cc_1$ so that if we choose $c_1$ such
that $ Cc_1=\frac14$ we have
$$ D_{\sigma}\subset K+\frac14B^{\infty}_{\sigma}$$
where $K:=C\epsilon^{-1}P_{\sigma}\Delta_{Nm}S.$  Now suppose $x\in
B^{\infty}_{\sigma}.$ Let
$a_1,a_2\in D_{\sigma} $ be defined by $a_1(j)=1$ if $x(j)\ge \frac12$
and
$a_1(j)=-1$ otherwise, while $a_2(j)=1$ if $x(j)\ge -\frac12$ and
$a_2(j)=-1$ otherwise.  Then
$$ \left\|x-\frac12(a_1+a_2)\r\|_{\infty}\le \frac12.$$
Thus
$$ B^{\infty}_{\sigma}\subset \Delta_2K+\frac34B_{\sigma}^{\infty}
= C\epsilon^{-1}P_{\sigma}\Delta_{2Nm}S +\frac34B_{\sigma}^{\infty}.$$
This implies for  $\theta = \frac34,$
$$ B^{\infty}_{\sigma} \subset
4 C\epsilon^{-1}\Gamma_{\theta} P_{\sigma}\Delta_{2Nm}S$$
 Letting $\varphi=\theta^{1/(2Nm)}$
and applying Lemma~\ref{approx2} we obtain
$$ \Gamma_{\theta}\Delta_{2Nm}S \subset \frac65 \Gamma_{\varphi}S.$$
 Note that $(\frac{3}{4})^{1/(2Nm)} \sim 1-(2Nm)^{-1} \ln (4/3) \le
1-c d^{-2}\epsilon^5(1-\ln\epsilon)^{-1}$ for some $c>0$ so that the
result follows.
\end{proof}

\begin{Thm}\label{pnormedmain}
 There is an absolute $C>0$ such that if
$\epsilon>0$ and  $X$ is a $p$-normed quasi-Banach space with $\dim X=n$
and $d(\hat X,\ell_{\infty}^n)\le d$ then $X$ has a quotient $Y$ with
$\dim Y\ge n(1-\epsilon)$ and $$d(Y,\ell_{\infty}^{\dim Y})\le
Cp^{-\frac{1}{p}}\epsilon^{4-\frac5p}(1-\ln \epsilon)^{\frac1p-1}
d^{\frac2p-1}.$$
\end{Thm}

\medskip
\noindent
{\bf Remark.} In \cite{K} examples are constructed of finite-dimensional
$p$-normed spaces $X_n$ (with $0<p<1$ fixed) so that $d(\hat
X_n,\ell_{\infty}^{\dim X_n})$ is uniformly bounded but
$\lim_{n\to\infty}\delta_{X_n}=\infty.$

\begin{proof} We can assume $B^{\infty}\subset B_{\hat X} \subset
dB^{\infty}.$  Then by Theorem \ref{main} we can find $\sigma$ with
$|\sigma|\ge n(1-\epsilon)$ so that
$$c\epsilon B^{\infty}_{\sigma}\subset \Gamma_{\theta}P_{\sigma}B_X$$
where $\theta= 1-cd^{-2}\epsilon^5(1-\ln \epsilon)^{-1}.$
Let $Y$ be the space of dimension $|\sigma|$ with unit ball
$B_Y=P_{\sigma}B_X.$  Since $B_Y$ is $p$-convex we have
(Lemma~\ref{pconv})
$$ \Gamma_{\theta}B_Y \subset p^{-\frac1p}
(cd^{-2}\epsilon^5(1-\log\epsilon)^{-1})^{1-\frac1p}B_Y.$$
 Finally observe that for a suitable $c>0$:
$$
cp^{\frac1p}d^{2-\frac2p}\epsilon^{\frac5p-4}(1-\log\epsilon)^{1-\frac1p}
B^{\infty}_{\sigma} \subset B_Y \subset d B^{\infty}_{\sigma}.$$
The result then follows.\end{proof}

\section{Cubic quotients}\label{cube2}

We start this section with the following lemma, which is in fact a corollary of
Theorem~\ref{main}.

\begin{Lem}\label{cubquo}
 Let $S$ be a compact spanning of $\mathbb R^n$  and $X$ be the
Banach space with unit ball $B_X = \Delta S.$ Let $m$ be
the largest integer such that $X$ has a subspace $Y$ of dimension
$m$ with $d(Y,\ell_1^m)\le 2$. Then for every integer $k$ satisfying
$2^{2k-1} \leq m$ there exists a rank $k$ projection $\pi$, so that for
some  cube $Q$ one has $Q \subset \Gamma_{b}\pi S \subset CQ$,
where $0<b<1$ is an absolute constant.
\end{Lem}

\begin{proof}
Let $Y$ be a subspace of $X$ of dimension $m$ so that $d (Y, \ell_1^m)\le 2.$
Then if $2^{2k-1}\le m$ there is a linear operator
$T:Y\to\ell_{\infty}^{2k}$ with $\|T\|\le 1$ and $T(B_Y)\supset \frac12
B^{\infty}_{2k}.$  $T$ can then be extended to a norm-one operator on $X$
and so $X$ has a quotient $Z$ of dimension $2k$ so that
$d(Z,\ell_{\infty}^{2k})\le 2.$   It follows immediately from
Theorem~\ref{main} with $\epsilon=\frac12$ that there is a further quotient $W$
of $Z$ with $\dim W \ge k$ and for some cube $Q_0$ in $W$, and fixed
constants $0<b<1$ and $1<C<\infty$, we have
$Q_0\subset \Gamma_{b}\pi_WS \subset CQ_0$ where $\pi_W$ is the quotient
map onto $W.$
\end{proof}

\begin{Thm}\label{final}  There is an absolute constant $c>0$ so that if
$X$ is a finite-dimensional $p$-normed space, then $X$ has a quotient
$E$ with $d(E,\ell_{\infty}^{\dim E})\le (cp)^{-1/p}$  and
$\dim E\ge c\ln A / (\ln \ln A )$, where $A= (p^{1/p}  \delta_X /4)^{p/(1-p)}$
(assuming that  $\delta_X$ is large enough).
\end{Thm}

\medskip

\noindent
{\bf Remark.} Take $X=\ell_p^n$ so that $\delta_X=n^{-1+1/p}.$  Then if
$X$ has a quotient $E$ of dimension $k$ with $d(E,\ell_{\infty}^k)\le C_p$
then $\hat X=\ell_1^n$ also has such a quotient which implies
$k\le c C_p \ln n = c C_p \ln \left(  \delta_X ^{p/(1-p)} \r) $. We conjecture
that this estimate is optimal up to an absolute constant, i.e. that every
$p$-normed space has a cubical quotient of such dimension. As one
can see from the proof below we could obtain such an estimate
(up to constant depending on $p$ only) if we were able to  prove the inclusion
with $\theta = 1- c (m \ln m) ^{-1}$  in Proposition~\ref{type2}.

\begin{proof} Let $S=B_X$ and  $m$ be the largest integer such that $X$
has a subspace $Y$ of dimension  $m$ with $d(Y,\ell_1^m)\le 2$.

Assume first $ m \leq 2^{2k}$. By Proposition~\ref{type2}
(and its proof) we have $ \Delta B_X \subset 4  \Gamma _{\theta} B_X$
for $\theta = 2^{- 1/N_k}$, where $N_k = (Ck) ^{C \ln \ln (Ck)}$. Then, by
Lemma~\ref{pconv}, we obtain
$$\Delta B_X \subset 4 p^{-1/p} ( 2 N_k)^{-1+1/p}$$
which implies
$$\delta _X \leq 4 p^{-1/p} ( 2 N_k)^{-1+1/p}.$$
Therefore $2 N_k \geq A := (p^{1/p} \delta _X /4 )^{p/(1-p)}$.
Finally we obtain
$k \geq C' \ln A / (\ln \ln A )$ (of course we may assume that $A > e^2$).

Suppose now $k \leq C' \ln A / (\ln \ln A )$. By above we have
$ m \geq 2^{2k}$. So Lemma~\ref{cubquo} implies the existence of
absolute constants $b$, $C_1$ and a rank $k$ projection $\pi$ such that
$Q\subset \Gamma_b \pi B_X \subset C_1 Q$ for some cube $Q$.
By Lemma~\ref{pconv} we obtain
$$
\Gamma_b \pi B_X \subset p^{-1/p} \left( 1- b \r) ^{1-1/p} \pi B_X
$$
so that  we have (if $E=X/\pi^{-1}(0)$),
$$ d(E,\ell_{\infty}^k) \le C_1 p^{-1/p}\left(1-b\r)^{1-1/p}.$$
This implies the theorem.
\end{proof}

\medskip

\noindent
{\bf Acknowledgment. }
The work on this paper was started during the visit of the second
named author to University of Missouri, Columbia.

\end{document}